\documentclass[11pt]{article} 

\newread\epsffilein    
\newif\ifepsffileok    
\newif\ifepsfbbfound   
\newif\ifepsfverbose   
\newif\ifepsfdraft     
\newdimen\epsfxsize    
\newdimen\epsfysize    
\newdimen\epsftsize    
\newdimen\epsfrsize    
\newdimen\epsftmp      
\newdimen\pspoints     
\def\epsfbox#1{\global\def\epsfllx{72}\global\def\epsflly{72}%
   \global\def\epsfurx{540}\global\def\epsfury{720}%
   \def\lbracket{[}\def\testit{#1}\ifx\testit\lbracket
   \let\next=\epsfgetlitbb\else\let\next=\epsfnormal\fi\next{#1}}%
\def\epsfgetlitbb#1#2 #3 #4 #5]#6{\epsfgrab #2 #3 #4 #5 .\\%
   \epsfsetgraph{#6}}%
\def\epsfnormal#1{\epsfgetbb{#1}\epsfsetgraph{#1}}%
\def\epsfgetbb#1{%
\openin\epsffilein=#1
\ifeof\epsffilein\errmessage{I couldn't open #1, will ignore it}\else
   {\epsffileoktrue \chardef\other=12
    \def\do##1{\catcode`##1=\other}\dospecials \catcode`\ =10
    \loop
       \read\epsffilein to \epsffileline
       \ifeof\epsffilein\epsffileokfalse\else
          \expandafter\epsfaux\epsffileline:. \\%
       \fi
   \ifepsffileok\repeat
   \ifepsfbbfound\else
    \ifepsfverbose\message{No bounding box comment in #1; using defaults}\fi\fi
   }\closein\epsffilein\fi}%
\def\epsfclipoff{\def\epsfclipstring{\ifepsfdraft\space clip\fi}}%
\epsfclipoff
\def\epsfsetgraph#1{%
   \epsfrsize=\epsfury\pspoints
   \advance\epsfrsize by-\epsflly\pspoints
   \epsftsize=\epsfurx\pspoints
   \advance\epsftsize by-\epsfllx\pspoints
   \epsfxsize\epsfsize\epsftsize\epsfrsize
   \ifnum\epsfxsize=0 \ifnum\epsfysize=0
      \epsfxsize=\epsftsize \epsfysize=\epsfrsize
      \epsfrsize=0pt
     \else\epsftmp=\epsftsize \divide\epsftmp\epsfrsize
       \epsfxsize=\epsfysize \multiply\epsfxsize\epsftmp
       \multiply\epsftmp\epsfrsize \advance\epsftsize-\epsftmp
       \epsftmp=\epsfysize
       \loop \advance\epsftsize\epsftsize \divide\epsftmp 2
       \ifnum\epsftmp>0
          \ifnum\epsftsize<\epsfrsize\else
             \advance\epsftsize-\epsfrsize \advance\epsfxsize\epsftmp \fi
       \repeat
       \epsfrsize=0pt
     \fi
   \else \ifnum\epsfysize=0
     \epsftmp=\epsfrsize \divide\epsftmp\epsftsize
     \epsfysize=\epsfxsize \multiply\epsfysize\epsftmp
     \multiply\epsftmp\epsftsize \advance\epsfrsize-\epsftmp
     \epsftmp=\epsfxsize
     \loop \advance\epsfrsize\epsfrsize \divide\epsftmp 2
     \ifnum\epsftmp>0
        \ifnum\epsfrsize<\epsftsize\else
           \advance\epsfrsize-\epsftsize \advance\epsfysize\epsftmp \fi
     \repeat
     \epsfrsize=0pt
    \else
     \epsfrsize=\epsfysize
    \fi
   \fi
   \ifepsfverbose\message{#1: width=\the\epsfxsize, height=\the\epsfysize}\fi
   \epsftmp=10\epsfxsize \divide\epsftmp\pspoints
   \vbox to\epsfysize{\vfil\hbox to\epsfxsize{%
      \ifnum\epsfrsize=0\relax
        \includegraphics{\ifepsfdraft}%
      \else
        \epsfrsize=10\epsfysize \divide\epsfrsize\pspoints
        \includegraphics{\ifepsfdraft}%
      \fi
      \hfil}}%
\global\epsfxsize=0pt\global\epsfysize=0pt}%
{\catcode`\%=12 \global\let\epsfpercent=
\long\def\epsfaux#1#2:#3\\{\ifx#1\epsfpercent
   \def\testit{#2}\ifx\testit\epsfbblit
      \epsfgrab #3 . . . \\%
      \epsffileokfalse
      \global\epsfbbfoundtrue
   \fi\else\ifx#1\par\else\epsffileokfalse\fi\fi}%
\def\epsfempty{}%
\def\epsfgrab #1 #2 #3 #4 #5\\{%
\global\def\epsfllx{#1}\ifx\epsfllx\epsfempty
      \epsfgrab #2 #3 #4 #5 .\\\else
   \global\def\epsflly{#2}%
   \global\def\epsfurx{#3}\global\def\epsfury{#4}\fi}%
\def\epsfsize#1#2{\epsfxsize}
\newtheorem{theorem}{Theorem}[section]
\newtheorem{lemma}[theorem]{Lemma}

\newtheorem{proposition}[theorem]{Proposition}

\newcommand{\lk}{\mbox{\rm Link}}
\newcommand{\tlk}{\mbox{\rm Tlk}}
\newcommand{\alk}{\mbox{\rm Alk}}

\newcommand{\glk}{\mbox{\rm Glk}}

\newcommand{\Int}{\mbox{\rm Int}}
\newcommand{\qed}{\hbox{\rule{6pt}{6pt}}}

\newcommand{\cmapright}[2]{
\smash{\mathop{
\hbox to 1cm{\rightarrowfill}}\limits^{#1}_{#2}}}

\newcommand{\cmapleft}[2]{
\smash{\mathop{
\hbox to 1cm{\leftarrowfill}}\limits^{#1}_{#2}}}

\begin{document}

\title{Triple Linking of Surfaces in $4$-Space{\footnote{MRCN:57Q45}}}

\author{
J. Scott Carter \\
University of South Alabama \\
Mobile, AL 36688 \\ carter@mathstat.usouthal.edu
 \and
Seiichi Kamada \\
Osaka City University \\
Osaka 558-8585, JAPAN\\ kamada@sci.osaka-cu.ac.jp \\
skamada@mathstat.usouthal.edu
\and
Masahico Saito \\
University of South Florida \\
Tampa, FL 33620 \\ saito@math.usf.edu
\and
Shin Satoh \\
Osaka City University \\
Osaka 558-8585, JAPAN\\
susato@sci.osaka-cu.ac.jp
}
\maketitle


\begin{abstract}
Triple linking numbers were defined for
$3$-component oriented surface-links in $4$-space using
signed triple points on projections in $3$-space.
In this paper we give an algebraic formulation using
intersections of homology classes (or cup products
on cohomology groups).  We prove that spherical links have
trivial triple linking numbers and that
triple linking numbers are link homology
invariants.
\end{abstract}


\section{Introduction}

A {\it surface-link} is
a closed surface $F$ embedded in ${\bf R}^4$
locally flatly.
In this paper,
we always assume that
$F$ is {\it oriented},
that is, each component of $F$ is
orientable and given a fixed orientation.
For a 3-component surface-link
$F=K_1\cup K_2\cup K_3$,
a linking number was defined in \cite{CJKLS} using its
projection in ${\bf R}^3$
in a way that is analogous to
the linking
number in classical knot theory.
In that paper,
it was introduced as
an example  of non-triviality of the state-sum invariants of
surface-links.
In fact,
the state-sum invariants in the classical link and
surface-link case
generalize linking number and Fox's coloring number.

In the current paper,
we give several alternative definitions of
the triple linking number and some properties.
The reader will find that this invariant is a quite natural
generalization of the notion of classical linking number
in contrast to a statement in
Rolfsen
\cite{Rolfsen} page 136:
``There is, however no analogous notion of linking number
to help us with codimension two link theory, for example,
in higher dimensions''.
We note, however, that Rolfsen himself with Massey
\cite{MassayRolf}
and with Fenn \cite{FennRolf} generalized classical linking
numbers to higher dimensions using degrees of maps.

Let $F=K_1\cup K_2\cup K_3$ be
a $3$-component surface-link in ${\bf R}^4$.
It is known (see \cite{CS:book} for example) that
a projection of $F$ into ${\bf R}^3$
can be assumed to have transverse double curves
and isolated branch/triple points. At a triple point, three
sheets intersect that have distinct relative heights with respect to the
projection direction, and we call them  top, middle, and bottom sheets,
accordingly.
If the orientation normals to the top, middle, bottom sheets at
a triple point
$\tau$ matches with this order the fixed orientation of ${\bf R}^3$,
then the sign of $\tau$ is positive and $\varepsilon(\tau)=1$.
Otherwise the sign is negative and $\varepsilon(\tau)=-1$.
 (See  \cite{CJKLS, CS:book}.) 
It is also known that any
closed 
oriented embedded surface $F$ in ${\bf R}^4$
bounds an oriented compact $3$-manifold $M$
embedded in ${\bf R}^4$, 
called
a {\it Seifert hypersurface} of $F$, such that $\partial M=F$.

We give six methods for defining an integer
(triple linking number) as follows.

\begin{itemize}
\item[(1)]
Consider a surface diagram of $K_1\cup K_2\cup K_3$
in ${\bf R}^3$.
A triple point is {\it of type $(i,j,k)$}
if the top sheet comes from $K_i$,
the middle comes from $K_j$, and
the bottom comes from $K_k$.
The sum of the signs
of all the triple points of type
$(1,2,3)$ is denoted by $\tlk_1(K_1,K_2,K_3)$.
This is the definition given in \cite{CJKLS}.

\item[(2)]
Let $M_i$ be a Seifert hypersurface for $K_i$ $(i=1,3)$.
Assume that $M_i\cap K_2$ is
a $1$-manifold in $K_2$ and that
$M_1\cap K_2$ and $M_3\cap K_2$
intersect transversely.
Count the intersections
between them algebraically
and denote the sum
by $\tlk_2(K_1,K_2,K_3)$.

\item[(3)]
Consider a Seifert hypersurface $M_1$ for $K_1$.
Assume that $M_1\cap K_2$ is a $1$-manifold,
which is disjoint from $K_3$.
The linking number $\lk(M_1\cap K_2,K_3)$
is denoted by $\tlk_3(K_1,K_2,K_3)$.

\item[(4)]
Let $M_i$ be a Seifert hypersurface for $K_i$ $(i=1,3)$
such that $M_1\cap M_3$ is a $2$-manifold
which intersects $K_2$ transversely.
Count the
intersections
between them algebraically
and denote
the sum
by $\tlk_4(K_1,K_2,K_3)$.

\item[(5)]
Let $M_i$ be a Seifert hypersurface for $K_i$
$(i=1,2,3)$ and let $N_2$ be a regular neighbourhood of $K_2$
in ${\bf R}^4$.
We may assume that $M_i\cap \partial N_2$
is a $2$-manifold in $\partial N_2$ and that
$M_1\cap \partial N_2$, $M_2\cap \partial N_2$ and
$M_3\cap \partial N_2$ intersect transversely
in a finite number of points.
Count the
intersections
algebraically
and denote the sum
by $\tlk_5(K_1,K_2,K_3)$.

\item[(6)]
Let $f:F_1 \cup F_2 \cup F_3 \rightarrow {\bf R}^4$ denote an
embedding of the disjoint union of oriented surfaces $F_i$
representing $F=K_1\cup K_2\cup K_3$.
Define a map
$L:F_1 \times F_2 \times F_3 \rightarrow S^3 \times S^3$ by
$$L(x_1,x_2,x_3) = \left( \frac {f(x_1)-f(x_2)}{||f(x_1)-f(x_2)||},
\frac {f(x_2)-f(x_3)}{||f(x_2)-f(x_3)||} \right)$$
for $x_1\in F_1$, $x_2\in F_2$ and $x_3 \in F_3$,
and denote the degree of $L$ by
$\tlk_6(K_1,K_2,K_3)$.

\end{itemize}


\begin{theorem}
$\tlk_i(K_1,K_2,K_3)=\pm\tlk_j(K_1,K_2,K_3)$
for any $i,j=1,\dots,6.$
\end{theorem}

\noindent
{\bf Remark.} \
In general,
the {\it triple linking number} $\tlk(K_i,K_j,K_k)$
for $i\ne j\ne k$ is defined to be the sum of
the signs of all the triple points of type $(i,j,k)$
on a surface diagram of $F$;

$$\tlk(K_i,K_j,K_k)=
\sum_{{\scriptstyle \tau:\
\mbox{\footnotesize type}\ (i,j,k)}}
\ \varepsilon(\tau).$$

\noindent
It is proved in \cite{CJKLS} that
this number is an invariant of the surface-link $F$
(independent of a diagram in ${\bf R}^3$)
by use of Roseman moves (Reidemeister moves for
surface-link diagrams) \cite{Rose},
and that this invariant vanishes in the case that
$i=k$;
that is, $\tlk(K_i,K_j,K_i)=0$
for $i\neq j$.
Hence throughout this paper,
we always assume that
$i, j, k$ are all distinct
whenever we refer to $\tlk(K_i,K_j,K_k)=\tlk_1(K_i,K_j,K_k)$.

\bigskip

We prove the following properties
of triple linking by using the above interpretations.


\begin{theorem}[\cite{CJKLS}]
\begin{itemize}
\setlength{\itemsep}{-3pt}
\item[{\rm (i)}]
$\tlk(K_1,K_2,K_3)=-\tlk(K_3,K_2,K_1)$.
\item[{\rm (ii)}]
$\tlk(K_1,K_2,K_3)+\tlk(K_2,K_3,K_1)+\tlk(K_3,K_1,K_2)=0$.
\end{itemize}
\end{theorem}

\begin{theorem}
\begin{itemize}
\setlength{\itemsep}{-3pt}
\item[{\rm (i)}]
If $K_2$ is homeomorphic to a $2$-sphere,
then $\tlk(K_1,K_2,K_3)=0$.
\item[{\rm (ii)}]
If both of $K_1$ and $K_3$ are homeomorphic to a $2$-sphere,
then $\tlk(K_1,K_2,K_3)=0$.
\end{itemize}
\end{theorem}

In \cite{NSato}
the asymmetric linking number $\alk(K,K')$ for
a two component oriented surface-link $F=K \cup K'$
was defined to be the non-negative generator
of the image of
$H_1(K) \rightarrow H_1(S^4 \backslash K') \cong {\bf Z}$.

\begin{theorem}
If $\alk(K_2, K_3)=0$,
then $\tlk(K_1, K_2,K_3)= \tlk(K_3, K_2, K_1)= 0$.
\end{theorem}

Two $n$-component surface-links
$F=K_1\cup \dots \cup K_n$ and $F'=K'_1\cup \dots \cup K'_n$
are {\it link homologous\/} if there is a
compact oriented $3$-manifold $W$ properly embedded in  ${\bf R}^4
\times [0,1]$ such that $W$ has $n$ components $W_1, \dots, W_n$
with $\partial W_i= K_i \times \{0\} \cup (-K'_i)\times
\{1\}$.
This relation is sometimes called {\it link-cobordism},
but that term also denotes the concordance relation.
Since link homotopy implies link homology, the following
theorem implies that triple linking invariants are
link homotopy invariants (this fact is also seen
from the sixth definition of $\tlk$).
For related topics, refer to
\cite{Tim1, TO, Kirk, KK, Koschorke,
MassayRolf, Ruberman, NSato}.


\begin{theorem}
Triple linking invariants are link homology invariants:
If $F=K_1\cup K_2\cup K_3$ and
$F'=K'_1\cup K'_2\cup K'_3$ are link homologous,
then $\tlk(K_i,K_j,K_k)=\tlk(K'_i,K'_j,K'_k)$.
\end{theorem}
See Remark~6.2 for further information about link homology.

By Theorem 1.2,
for any 3-component surface-link
$F=K_1\cup K_2\cup K_3$,
there exists a pair of integers $a$ and $b$ such that

$$(*) \left\{
\begin{array}{l}
\tlk(K_1,K_2,K_3)=-\tlk(K_3,K_2,K_1)=-(a+b),\\
\tlk(K_2,K_3,K_1)=-\tlk(K_1,K_3,K_2)=b, \\
\tlk(K_3,K_1,K_2)=-\tlk(K_2,K_1,K_3)=a.
\end{array} \right.$$

\noindent
In \cite{CJKLS}, it is shown that for any pair of
integers $a$ and $b$, there exists a surface-link $F$
whose triple linking numbers satisfy the above equations.
However, that paper does not treat any problem
about genera of the components of $F$.
By Theorem 1.3, we see that

(1) if $a\ne 0$ and $b=0$, then
$g(K_i)\geq 1$ $(i=1,2)$, and

(2) if $a\ne 0$, $b\ne 0$ and $a+b\ne 0$,
then $g(K_i)\geq 1$ $(i=1,2,3)$,

\noindent
where $g(K_i)$ denotes the genus of $K_i$.

\begin{proposition}
{\rm (i)} For any integer $a\ne 0$,
there exists a surface-link
$F=K_1\cup K_2\cup K_3$ whose
triple linking numbers satisfy the above equations $(*)$
with $b=0$
and $g(K_i)=1$ $(i=1,2)$ and $g(K_3)=0$.

{\rm (ii)}
For any pair of integers $a$ and $b$
with $a\ne 0$, $b\ne 0$ and $a+b\ne 0$,
there exists a surface-link
$F=K_1\cup K_2\cup K_3$ whose
triple linking numbers satisfy the above equations $(*)$
and $g(K_i)=1$ $(i=1,2,3)$.
\end{proposition}

This paper is organized as follows:
in Section 2,
we interprete $\tlk_1$ in terms of
the decker curves of a surface diagram.
In Section 3 we give precise definitions of
triple linking numbers $\tlk_i$ for $i=2,\dots,5$
(in terms of homology)
and prove Theorem 1.1.
Section 4 is devoted to
proving Theorems 1.2--1.5.
Proposition 1.6 is proved in Section 5.

\bigskip

Throughout this paper,
all the homology and cohomology groups
have the {\bf Z}-coefficient.


\section{Decker Curves and Triple Linking}

Let $F$ be a surface-link
and $F^*$ a surface diagram of $F$
with respect to a projection
$p:{\bf R}^4\rightarrow{\bf R}^3$.
Let $\Gamma(F^*)$
denote the double point set of $F^*$;
$$\{ p(x)\ |\  x \in F, \ p(x)=p(y) \ \mbox{for some} \
  y \in F, x \neq y \},$$
which consists of immersed curves,
called {\it double curves}.
A double curve $C^*$ is an immersed circle
or an immersed arc in ${\bf R}^3$.
If $C^*$ is an immersed circle, then
$(p|_F)^{-1}(C^*)=C\cup C'$ for
some pair of immersed circles $C$ and $C'$ in $F$.
If $C^*$ is an immersed arc,
then its endpoints are branch points of $F^*$ and
$(p|_F)^{-1}(C^*)=C\cup C'$ for
some pair of immersed arcs $C$ and $C'$ in $F$
with $\partial C= \partial C'$.
The curves $C$ and $C'$ are called {\it decker curves\/}
over $C^*$:
one of them is in higher position
than the other with respect to the projection direction,
which is called an {\it upper decker curve\/}
 and the other
is called
an {\it lower decker curve\/}.
We notice that
the preimage of a triple point
consists of three points of $F$
which are intersections of decker curves.
See \cite{CS:book} for details.

Double curves and decker curves are oriented as
follows:
Let $x$ be a point of $F$
whose image $x^*=p(x)$ is not a branch point.
There is a regular neighborhood
$N$ of $x$ in $F$ such that $p|_N$ is an embedding.
An {\it orientation normal
$\vec{n}$ to $N^*=p(N)$ \/}
in ${\bf R}^3$ at $x^*$
is specified in such a way that
$(\vec{v}_1, \vec{v}_2, \vec{n})$
matches the orientation of ${\bf R}^3$,
where the pair of tangents
$(\vec{v}_1, \vec{v}_2)$ defines
the orientation of $N^*$ that is
induced from the orientation of $N \subset F$.
If $y$ is a double point on a double curve $C^*$, then
$C^*$ is locally an intersection of $N_1^*$ and
$N_2^*$, where $N_1^*$ is upper and $N_2^*$ is lower.
We assign a tangent vector $\vec{v}$ of $C^*$
at $y$ such that
$(\vec{n}_1,
\vec{n}_2, \vec{v})$ matches the orientation of ${\bf R}^3$.
This defines an orientation of $C^*$, cf.
\cite{CJKLS, CS:book}. 
We give an orientation
to the lower decker
 curve over $C^*$ such that it inherits
the orientation from $C^*$,
and give the opposite orientation to
the upper decker curve.
Note that, if $C^*$ is an arc, then the orientations of
$C$ and $C'$ are compatible (i.e., the union $C \cup C'$
forms an oriented immersed circle in $F$).

Let $F=K_1 \cup K_2 \cup K_3$
be a 3-component surface-link.
A double curve $C^*$ is
{\it of type $(i,j)$} if the upper decker
curve lies in $K_i$ and the lower decker
curve lies in $K_j$.
A decker curve over $C^*$ is {\it of type $(i,j)$} if
$C^*$ is so.

At a triple point $\tau$,
if the orientation normals to the top, middle,
and bottom sheets at $\tau$ matches with this order
the fixed orientation of ${\bf R}^3$,
then the sign of $\tau$ is positive and
$\varepsilon(\tau)=+1$;
otherwise the sign is negative and $\varepsilon(\tau)=-1$.

We interprete the triple linking $\tlk_1$
in terms of double decker curves as follows.
Let $D_{12}^{\ell}$ (resp. $D_{23}^u$)
denote the union of
lower decker curves of type $(1,2)$ (resp.
upper decker curves of type $(2,3)$).
Note that both $D_{12}^{\ell}$ and $D_{23}^u$
are contained in $K_2$.

\begin{lemma}
$\tlk_1(K_1,K_2,K_3)=-\Int_{K_2}(D_{12}^{\ell},D_{23}^u),
$
where $\Int_{K_2}(D_{12}^{\ell},D_{23}^u)$
is the intersection number in $K_2$.
\end{lemma}

\noindent
{\it Proof.}
Let $\tau$ be a triple point
of type $(1,2,3)$.
The preimage of $\tau$ consists of three points of $F$.
Exactly one of them is
on $K_2$ and that is a double point
of $D_{12}^{\ell}$ and $D_{23}^u$.
Conversely the image of
a double point of  $D_{12}^{\ell}$ and $D_{23}^u$ is
a triple point of $F^*$ of type $(1,2,3)$.
Hence there is a one-to-one correspondence
between the set of triple points of type $(1,2,3)$ and
double points of $D_{12}^{\ell}$ and $D_{23}^u$.
If the sign of $\tau$ is positive (or negative, resp.)
then the corresponding intersection of $D_{12}^{\ell}$ and $D_{23}^u$
is negative (resp. positive),
see Figure~1.  Thus we have the result.
\ $\qed$




\begin{figure}[htb]
\begin{center}
\mbox{
\epsfxsize=5in
\epsfbox{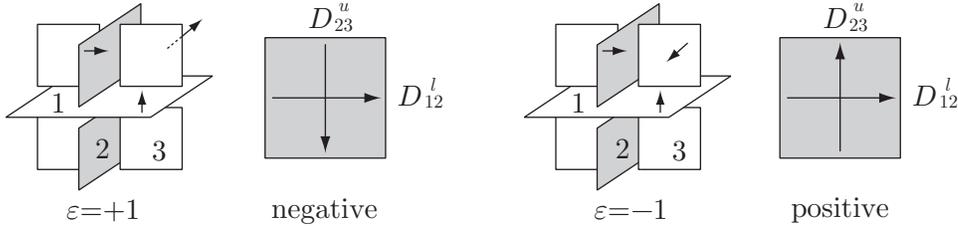}
}
\end{center}
\caption{Triple points and intersection of decker curves}
\label{Fig01}
\end{figure}


Since $D_{12}^{\ell}$ is the union of circles in ${\bf R}^4$
disjoint from $K_3$,
the linking number $\lk(D_{12}^{\ell}, K_3)$
is defined.


\begin{lemma}
$\tlk_1(K_1,K_2,K_3)=\lk(D_{12}^{\ell}, K_3).$
\end{lemma}

\noindent
{\it Proof.}
Without loss of generality, we may assume that
the projection $p$ is given by
$p(w,x,y,z) = (x,y,z)$.
For a real number $\lambda$ we denote by
$t_\lambda: {\bf R}^4 \to {\bf R}^4$ the translation
with $t_\lambda(w,x,y,z)= (w+\lambda, x,y,z)$.
Let $M'_3$ be a 3-chain in ${\bf R}^4$
with $\partial M'_3 = K_3$.
Take a sufficiently large number $R$ and
consider a $3$-chain

$$M_3 = \cup_{\lambda \in [0, R]} t_\lambda(K_3)
+ t_{R}(M'_3)   $$

\noindent
so that $\partial M_3= K_3$ and
$D_{12}^{\ell} \cap M_3 = D_{12}^{\ell}
\cap (\cup_{\lambda \in [0, R]}
t_\lambda(K_3))$.
The projection $p$ induces a one-to-one correspondence between
the geometric intersection
$D_{12}^{\ell} \cap (\cup_{\lambda \in [0, R]}
t_\lambda(K_3))$ and the subset of
$D_{12}^{\ell *} \cap K_3^* = p(D_{12}^{\ell}) \cap p(K_3)$
consisting of points where $D_{12}^{l*}$ is higher than $K_3^*$
(in the over-under information of the surface diagram $F^*$),
i.e., the set of triple points of $F^*$ of type $(1,2,3)$.
Since the orientation of  $D_{12}^{\ell}$ is parallel
to the orientation of $D_{12}^{\ell *}$, the sign of an
intersection of $D_{12}^{\ell}$ and $\cup_{\lambda \in [0, R]}
t_\lambda(K_3)$ coincides with the sign of
the corresponding intersection of $D_{12}^{\ell*}$
and $K_3^*$, which is the sign of the triple point
(see Figure~2).
Thus we have the result.
\ \qed




\begin{figure}[htb]
\begin{center}
\mbox{
\epsfxsize=5in
\epsfbox{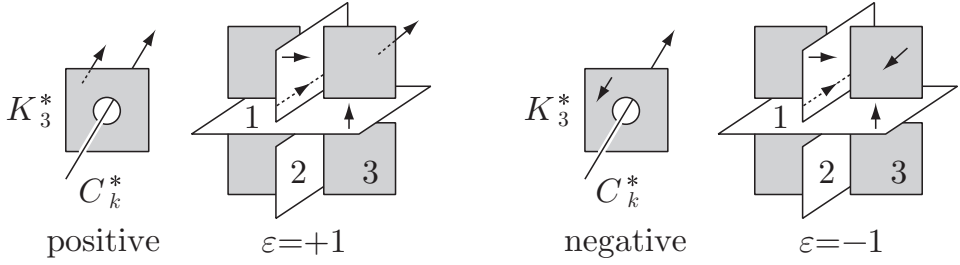}
}
\end{center}
\caption{The intersection between $K_3$ and $D^\ell_{12}$}
\label{Fig02}
\end{figure}



\section{Proof of Theorem 1.1}

For a compact oriented $n$-manifold $M$ with
$\{A,B\}= \{\partial M, \emptyset\}$,
we denote by
$$
\cdot_M \quad : H_p(M,A) \times H_q(M,A) \rightarrow
H_{p+q-n}(M,A)
$$
the intersection map, which is defined by
$$
x \cdot_M y = P_M( P_M^{-1}(x) \cup P_M^{-1}(y) )
$$
where $P_M: H^{\ast}(M,B) \rightarrow H_{n-\ast}(M,A)$ is
the Poincar\'e duality isomorphism
(see \cite{Massey}, page~391).
We will use $\cdot$ and $P$ instead of $\cdot_M$ and $P_M$
when their meanings are obvious in context.

Let $F=K_1\cup K_2\cup K_3$
be a $3$-component surface-link.
For simplicity of argument, we
assume that $F$ is embedded in
the $4$-sphere
$S^4={\bf R}^4\cup\{\infty\}$.
For a regular neighborhood $N_i$ of $K_i$
in $S^4$,
we put
$$E_i={\rm Cl}(S^4 \backslash N_i), \
E_{ij}={\rm Cl}(S^4 \backslash (N_i\cup N_j))
\mbox{ for $i \ne j$, and }
E= {\rm Cl}(S^4 \backslash (N_1\cup N_2\cup N_3)),$$
where Cl denotes the closure.
We denote by $M_i$ a $3$-chain in $S^4$ with $\partial M_i = K_i$
for $i=1,2,3$ (the reader may suppose that
it is a Seifert hypersurface
for $K_i$,  i.e., a compact oriented $3$-manifold
embedded in $S^4$ with $\partial M_i=K_i$).
We also denote by $M_i$ the homology class in
$$H_3(S^4,K_i) {\cong} H_3(S^4,N_i) {\cong}
H_3(E_i,\partial E_i)$$
represented by $M_i$.
By $u_i \in H^1(E_i)$ we denote the Poincar\'e dual of
$M_i \in H_3(E_i,\partial E_i)$, i.e.,
$M_i= P(u_i) = u_i\cap[E_i]$.
For a subset $X$ of $E_i$, we will denote by
$u_i|_{X} \in H^1(X)$ the image of $u_i$ by the
inclusion-induced homomorphism
$H^1(E_i) \to H^1(X)$.
Moreover, if $X$ is an $n$-manifold, we denote by
$M_i|_{(X, \partial X)}$
(or $M_i|_{X}$ if $\partial X =\emptyset$)
the Poincar\'e dual $P_X(u_i|_{X})= (u_i|_{X}) \cap [X] \in
H_{n-1}(X, \partial X)$ of $u_i|X$.

\bigskip

For $i \in \{1,3\}$, since $K_2 \subset E_i$,
$M_i|_{K_2} \in H_1(K_2)$  is defined. (When
we consider $M_i$ as
a 
$3$-chain,
the intersection of $M_i$ and $K_2$
(as a $1$-cycle in $K_2$)
represents
$M_i|_{K_2}$.)
Let
$$\tlk_2(K_1, K_2, K_3)=
\varepsilon_{K_2}(M_1|_{K_2} \cdot M_3|_{K_2}),$$
where $\varepsilon_{K_2}: H_0(K_2) \to {\bf Z}$
is the augmentation.


\begin{lemma}
$\tlk_1(K_1,K_2,K_3)=\tlk_2(K_1,K_2,K_3).$
\end{lemma}

\noindent
{\it Proof.}
We assume $F \subset {\bf R}^4 \subset S^4$ and continue
the situation of the proof of Lemma~2.2.
Let $M'_1$ be a 3-chain in ${\bf R}^4 (\subset S^4)$
with $\partial M'_1= K_1$ and
consider a $3$-chain $M_1$ such that
$$M_1 = -\cup_{\lambda \in [-R, 0]} t_\lambda(K_1)
+ t_{-R}(M'_1)$$
with $\partial M_1= K_1$.
The intersection of $M_1$ and $K_2$ is equal to
that of
$-\cup_{\lambda \in [-R, 0]} t_\lambda(K_1)$ and $K_2$ which
is the 1-chain $-D_{12}^{\ell}$ in $K_2$, and
the intersection of
$M_3$ and $K_2$ is equal to that of
$\cup_{\lambda \in [0, R]} t_\lambda(K_3)$ and
$K_2$  which is the 1-chain $D_{23}^u$ in $K_2$.
Therefore, by Lemma~2.1, we have

$$\begin{array}{ll}
\tlk_1(K_1,K_2,K_3)
&= -\Int_{K_2}(D_{12}^{\ell},D_{23}^u) \\
&= -\varepsilon_{K_2}(D_{12}^{\ell} \cdot D_{23}^u) \\
&= -\varepsilon_{K_2}(-M_1|_{K_2} \cdot M_3|_{K_2})\\
&= \tlk_2(K_1,K_2,K_3). \ \qed
\end{array}$$


\bigskip

\noindent
{\bf  Remark.} \
The argument in Lemmas~2.2 and 3.1 implies that
for any Seifert hypersurface $M_1$ for $K_1$,
the intersection of $M_1$ and $K_2$ (as a
1-cycle in $K_2$)
is homologous to $-D_{12}^{\ell}$,
 and that
for any Seifert hypersurface $M_3$ for $K_3$,
the intersection of $M_3$ and $K_2$ (as a 1-cycle in $K_2$)
is homologous to $D_{23}^u$.

\bigskip

We denote by $[M_1 \cap K_2]_{E_3} \in H_1(E_3)$ the
homology class of the intersection $M_1 \cap K_2$ as a $1$-cycle
in $E_3$. This is equal to the image of $M_1|_{K_2} \in H_1(K_2)$
under the inclusion-induced
homomorphism $H_1(K_2) \to H_1(E_3)$ and  also
equal to the image of
$M_1|_{(E_{13}, \partial E_{13})}
\cdot  [K_2]_{E_{13}} \in H_1(E_{13})$ under the
inclusion-induced
homomorphism $H_1(E_{13}) \to H_1(E_3)$,
where $[K_2]_{E_{13}} \in H_2(E_{13})$ is represented by $K_2$.
Let
$$ \begin{array}{ll}
\tlk_3(K_1, K_2, K_3)
&=
\lk([M_1 \cap K_2]_{E_3}, K_3) \\
&=
\varepsilon_{E_3} ([M_1 \cap K_2]_{E_3} \cdot M_3),
\end{array} $$
where $M_3 \in H_3(E_3, \partial E_3)$ and
$\varepsilon_{E_3}: H_0(E_3) \to {\bf Z}$ is the augmentation.


\begin{lemma}
$\tlk_1(K_1,K_2,K_3)=-\tlk_3(K_1,K_2,K_3).$
\end{lemma}

\noindent
{\it Proof.}
In the situation of the proof of Lemma~2.2,
$[M_1 \cap K_2]_{E_3} \in H_1(E_3)$ is
represented by the $1$-cycle $-D_{12}^{\ell}$.
Therefore, by Lemma~2.2, we have
$$\begin{array}{ll}
\tlk_3(K_1,K_2,K_3)
&= \lk([M_1 \cap K_2]_{E_3}, K_3) \\
&= \lk(-D_{12}^{\ell}, K_3) \\
&= -\tlk_1(K_1,K_2,K_3). \ \qed
\end{array}$$

\bigskip

We denote by $[M_1 \cap M_3]_{(E_{13}, \partial E_{13})}
\in H_1(E_{13}, \partial E_{13})$
the class  of the intersection $M_1 \cap M_2$ as a $2$-cycle in
$(E_{13}, \partial E_{13})$ when we regard $M_i$ as a $3$-chain.
This is equal to the intersection product
$M_1 |_{(E_{13}, \partial E_{13})}
\cdot M_3 |_{(E_{13}, \partial E_{13})}
\in H_2(E_{13}, \partial E_{13})$.
Let
$$\begin{array}{ll}
\tlk_4(K_1, K_2, K_3)
&=
\varepsilon_{E_{13}} ([M_1 \cap M_3]_{(E_{13}, \partial E_{13})}
  \cdot [K_2]_{E_{13}}) \\
&=
\varepsilon_{E_{13}} (M_1 |_{(E_{13}, \partial E_{13})}
\cdot M_3 |_{(E_{13}, \partial E_{13})} \cdot [K_2]_{E_{13}}).
\end{array}$$


\begin{lemma}
$\tlk_3(K_1,K_2,K_3)=\tlk_4(K_1,K_2,K_3).$
\end{lemma}

\noindent
{\it Proof.}
Let $i_*: H_*(E_{13}) \to H_*(E_3)$ and
$i^*: H^*(E_3) \to H^*(E_{13})$ be the inclusion-induced
homomorphisms.
Recall that
$[M_1 \cap K_2]_{E_3}= i_*(M_1|_{(E_{13}, \partial E_{13})}
\cdot  [K_2]_{E_{13}})$.  Thus,

$$\begin{array}{ll}
\tlk_4(K_1, K_2, K_3)
&=
\varepsilon_{E_{13}} (
M_1 |_{(E_{13}, \partial E_{13})}
\cdot M_3 |_{(E_{13}, \partial E_{13})}
\cdot [K_2]_{E_{13}}) \\
&=
\varepsilon_{E_{3}} \circ i_* (
M_1 |_{(E_{13}, \partial E_{13})}
\cdot M_3 |_{(E_{13}, \partial E_{13})}
\cdot [K_2]_{E_{13}})   \\
&= -
\varepsilon_{E_{3}} \circ i_* (
M_3 |_{(E_{13}, \partial E_{13})}
\cdot M_1 |_{(E_{13}, \partial E_{13})}
\cdot [K_2]_{E_{13}})  \\
&= -
\varepsilon_{E_{3}} \circ i_* (
u_3 |_{E_{13}}
\cap (M_1 |_{(E_{13}, \partial E_{13})}
\cdot [K_2]_{E_{13}}) )  \\
&= -
\varepsilon_{E_{3}} \circ i_* (
i^*(u_3)
\cap (M_1 |_{(E_{13}, \partial E_{13})}
\cdot [K_2]_{E_{13}}) )  \\
&= -
\varepsilon_{E_{3}}  (
u_3
\cap i_* (M_1 |_{(E_{13}, \partial E_{13})}
\cdot [K_2]_{E_{13}}) )  \\
&= -
\varepsilon_{E_{3}}  (
M_3
\cdot [M_1 \cap K_2]_{E_3} )  \\
&=
\varepsilon_{E_{3}}  (
[M_1 \cap K_2]_{E_3} \cdot M_3
)  \\
&=
\tlk_3(K_1, K_2, K_3). \ \qed
 \end{array}$$

\bigskip

For $i \in \{1,2,3\}$, since $\partial N_2 \subset E_i$,
$M_i|_{\partial N_2} \in H_2(\partial N_2)$ is defined.  Let

$$\begin{array}{ll}
\tlk_5(K_1,K_2,K_3)
&=
\varepsilon_{\partial N_2}(
M_1|_{\partial N_2} \cdot
M_2|_{\partial N_2} \cdot
M_3|_{\partial N_2}) \\
&=
<u_1|_{\partial N_2} \cup
u_2|_{\partial N_2} \cup
u_3|_{\partial N_2}, [\partial N_2] >.
\end{array}  $$


\begin{lemma}
$\tlk_4(K_1,K_2,K_3)=\tlk_5(K_1,K_2,K_3).$
\end{lemma}

\noindent
{\it Proof.}
Let $i : \partial N_2 \to N_2$ be the inclusion map.
In $H_0(N_2)$, we have
$$
\begin{array}{ll}
i_*(
M_1|_{\partial N_2} \cdot
M_2|_{\partial N_2} \cdot
M_3|_{\partial N_2})
&= -
i_*(
M_1|_{\partial N_2} \cdot
M_3|_{\partial N_2} \cdot
M_2|_{\partial N_2}) \\
&= -
i_*(
\partial_*(M_1|_{(N_2,\partial N_2)}) \cdot
\partial_*(M_3|_{(N_2,\partial N_2)}) \cdot
M_2|_{\partial N_2}) \\
&= -
i_*(
\partial_*(M_1|_{(N_2,\partial N_2)} \cdot
M_3|_{(N_2,\partial N_2)}) \cdot
M_2|_{\partial N_2}) \\
&= -
(M_1|_{(N_2,\partial N_2)} \cdot
M_3|_{(N_2,\partial N_2)}) \cdot
i_*(M_2|_{\partial N_2}) \\
&= -
M_1|_{(N_2,\partial N_2)} \cdot
M_3|_{(N_2,\partial N_2)} \cdot
(-[K_2]_{N_2}) \\
&=
M_1|_{(N_2,\partial N_2)} \cdot
M_3|_{(N_2,\partial N_2)} \cdot
[K_2]_{N_2}.
\end{array}
$$
Thus
$$
\tlk_5(K_1,K_2,K_3) = \varepsilon_{N_2}(
M_1|_{(N_2,\partial N_2)} \cdot
M_3|_{(N_2,\partial N_2)} \cdot
[K_2]_{N_2}
). $$
It is obvious that
$$
\varepsilon_{N_2}(
M_1|_{(N_2,\partial N_2)} \cdot
M_3|_{(N_2,\partial N_2)} \cdot
[K_2]_{N_2} )
=
\varepsilon_{E_{13}}(
M_1|_{(E_{13},\partial E_{13})} \cdot
M_3|_{(E_{13},\partial E_{13})} \cdot
[K_2]_{E_{13}} )
$$
and hence we have the result.
\ \qed

\bigskip


\begin{lemma}
${\rm Tlk}_6(K_1,K_2,K_3) =
\pm {\rm Tlk}_1(K_1,K_2,K_3).$
 \end{lemma}

\noindent
{\it Proof.\/}
Since ${\rm Tlk}_6$ is an ambient isotopy invariant,
we may assume that the surface-link
$F= f(F_1) \cup f(F_2) \cup
f(F_3)$  is in general position with respect to
the projection $p : {\bf R}^4 \to {\bf R}^3$
with  $p(w,x,y,z) = (x,y,z)$.
  The preimage of a particular point $((1,0,0,0), (1,0,0,0))$
by $L$ consists of triples
$(x_1,x_2,x_3) \in F_1 \times F_2 \times
F_3$  such that $p(f(x_1))= p(f(x_2)) = p(f(x_3))$ and
$f(x_1)$ is the upper, $f(x_2)$ is the middle, $f(x_3)$ is the lower
lift of the triple point $p(f(x_1))$.
   For each such triple $(x_1, x_2, x_3)$,  let
$D_T^2$, $D_M^2$, $D_B^2$ be regular neighborhoods of them in
$F_1 \cup F_2 \cup F_3$, and let $\varepsilon \in \{+1, -1\}$
be the sign of the triple point $p(f(x_1))$.
Let $(x_T,y_T)$, $(x_M, -\varepsilon z_M)$ and
$(y_B,z_B)$ be coordinate systems
of
$D_T^2$, $D_M^2$ and $D_B^2$ around $x_1, x_2$ and $x_3$,
respectively.
Modifying $f$ up to ambient isotopy, we may assume that
the restriction of $f$ to
$D_T^2 \cup D_M^2 \cup D_B^2$ is given by
defined by
\begin{eqnarray*}
(x_T,y_T) &
\mapsto & (0, x_0, y_0, z_0)+ (3,x_T,y_T,0) \\
(x_M,- \varepsilon z_M) &
\mapsto & (0, x_0, y_0, z_0)+ (2,x_M, 0, z_M) \\
(y_B,z_B) &
\mapsto & (0, x_0, y_0, z_0)+ (1,y_B, z_B) \end{eqnarray*}
where
$(x_0, y_0, z_0) \in {\bf R}^3$ is the triple point $p(f(x_1))$.
In this situation, the restriction
$$L': D_T^2 \times D_M^2 \times D_B^2 \rightarrow S^3 \times
S^3$$ is given by the formula
$$\left(
\frac
{(1,x_T-x_M,y_T,-z_M)}
{\sqrt{1+ (x_T-x_M)^2+ y_T^2 + z_M^2}},
\frac
{(1,x_M,-y_B,z_M-z_B)}
{\sqrt{1+ x_M^2 +y_B^2 +(z_M-z_B)^2}}
\right).$$
The map $L'$ is injective and hence it is
a homeomorphism onto its image.
Its (local) degree is $+1$ or $-1$ which depends only on
$\varepsilon$.
Since the degree of $L$ is the sum of the (local) degrees  of
$L'$ for all triples $(x_1, x_2, x_3)$ in the preimage
$L^{-1}((1,0,0,0),(1,0,0,0))$, this number
agrees up to sign with the
triple linking number ${\rm Tlk}_1(f(F_1), f(F_2), f(F_3)).$
\ \qed

\bigskip

By Lemmas 3.1--3.5,
we have Theorem 1.1.



\section{Proof of Theorems 1.2--1.5}

To prove Theorem 1.2,
it is useful to change $\partial N_2$
in the definition of $\tlk_5$
for $\partial E_2$.

\begin{lemma}
$\tlk_5(K_1,K_2,K_3)
=- \varepsilon_{\partial E_2}(
M_1|_{\partial E_2} \cdot
M_2|_{\partial E_2} \cdot
M_3|_{\partial E_2})$, where
the intersections
are taken in $\partial E_2$.
\end{lemma}

\noindent
{\it Proof.}
Since $\partial N_2$ and $\partial E_2$ are the
same $3$-submanifold of $S^4$ with opposite
orientations,
$[\partial N_2] = - [\partial E_2]$
in $H_3(\partial N_2)= H_3(\partial E_2)$.
Thus, in $H_0(\partial N_2)=H_0(\partial E_2)$,
$$
\begin{array}{ll}
M_1|_{\partial N_2} \cdot
M_2|_{\partial N_2} \cdot
M_3|_{\partial N_2}
&=
(u_1|_{\partial N_2} \cup
u_2|_{\partial N_2} \cup
u_3|_{\partial N_2}) \cup [\partial N_2]  \\
&=
(u_1|_{\partial E_2} \cup
u_2|_{\partial E_2} \cup
u_3|_{\partial E_2}) \cup (-[\partial E_2])  \\
&= -
M_1|_{\partial E_2} \cdot M_2|_{\partial E_2} \cdot
M_3|_{\partial E_2}.   \ \qed
\end{array}
$$

\noindent
{\it Proof of Theorem $1.2$} (i)
$$
\begin{array}{ll}
\tlk_2(K_1, K_2, K_3)
&=
\varepsilon_{K_2}(M_1|_{K_2} \cdot M_3|_{K_2}) \\
&= -
\varepsilon_{K_2}(M_3|_{K_2} \cdot M_1|_{K_2}) \\
&= -
\tlk_2(K_3, K_2, K_1).
\end{array} $$

(ii)
Note that
$M_i|_{(E, \partial E)} = P_{E}(u_i|_{E})
\in H_3(E, \partial E)$
is the image of $M_i$
under
$$
H_3(E_i, \partial E_i) \to
H_3(E_i, \partial E_i \cup N_j \cup N_k)
\cong
H_3(E, \partial E),$$  and
$M_i|_{\partial E_2}
\in H_2(\partial E_2)$ is
the image of $M_i|_{(E, \partial E)}$ under
$$H_3(E, \partial E) \to H_2(\partial E)
\cong
H_2(\partial E_1) \oplus H_2(\partial E_2) \oplus
H_2(\partial E_3) \to H_2(\partial E_2), $$
the boundary operator followed by the projection to
$H_2(\partial E_2)$.
    We denote by
$(M_i|_{\partial E_2})_{\partial E} \in H_2(\partial E)$
the image of $M_i|_{\partial E_2}$ under
the inclusion-induced homomorphism
$H_2(\partial E_2) \to H_2(\partial E)$.
By Lemma 4.1,
$$
\tlk_5(K_1, K_2, K_3)= -
\varepsilon_{\partial E}(
(M_1|_{\partial E_2})_{\partial E} \cdot
(M_2|_{\partial E_2})_{\partial E} \cdot
(M_3|_{\partial E_2})_{\partial E}).    $$
Thus, we have
$$\begin{array} {rl}
\tlk_5(K_1, K_2, K_3) + &
\tlk_5(K_2, K_3, K_1) +
\tlk_5(K_3, K_1, K_2) \\
=
& - \varepsilon_{\partial E}(
(M_1|_{\partial E_2})_{\partial E} \cdot
(M_2|_{\partial E_2})_{\partial E} \cdot
(M_3|_{\partial E_2})_{\partial E} \\
& +
(M_2|_{\partial E_3})_{\partial E} \cdot
(M_3|_{\partial E_3})_{\partial E} \cdot
(M_1|_{\partial E_3})_{\partial E} \\
& +
(M_3|_{\partial E_1})_{\partial E} \cdot
(M_1|_{\partial E_1})_{\partial E} \cdot
(M_2|_{\partial E_1})_{\partial E}
) \\
=
& - \varepsilon_{\partial E}(
(M_1|_{\partial E_2})_{\partial E} \cdot
(M_2|_{\partial E_2})_{\partial E} \cdot
(M_3|_{\partial E_2})_{\partial E} \\
& +
(M_1|_{\partial E_3})_{\partial E} \cdot
(M_2|_{\partial E_3})_{\partial E} \cdot
(M_3|_{\partial E_3})_{\partial E} \\
& +
(M_1|_{\partial E_1})_{\partial E} \cdot
(M_2|_{\partial E_1})_{\partial E} \cdot
(M_3|_{\partial E_1})_{\partial E}
) \\
=
& - \varepsilon_{\partial E}(
M_1|_{\partial E} \cdot
M_2|_{\partial E} \cdot
M_3|_{\partial E}) \\
=
& - \varepsilon_{\partial E}(
\partial_*(M_1|_{(E,\partial E)}) \cdot
\partial_*(M_2|_{(E,\partial E)}) \cdot
\partial_*(M_3|_{(E,\partial E)})) \\
=
& - \varepsilon_{\partial E}( \partial_* (
M_1|_{(E,\partial E)} \cdot
M_2|_{(E,\partial E)} \cdot
M_3|_{(E,\partial E)})) \\
=
& 0. \  \qed
\end{array}$$

\noindent
{\it Proof of Theorem $1.3$} \
(i) The intersection number
between two oriented curves on a 2-sphere vanishes.
By Lemma 2.1, we have
$\tlk(K_1,K_2,K_3)=0$.

(ii) This is an immediate consequence of
(i) and Theorem 1.2(ii).
\ $\qed$

\bigskip

\noindent
{\it Proof of Theorem $1.4$} \
If $\alk(K_2, K_3)=0$, then
$\tlk_3(K_1, K_2, K_3) =\lk([M_1 \cap K_2]_{E_3}, K_3) = 0$,
for $[M_1 \cap K_2]_{E_3} \in H_1(E_3)$  is the image of
$M_1|_{K_2} \in H_1(K_2)=0$.
By Theorem~1.2, we have
$\tlk(K_3, K_2, K_1)=0$.
\ $\qed$

\bigskip

We consider surface-links $F=K_1\cup K_2\cup K_3$
in which each $K_i$ is not necessarily connected.
Such a surface-link is called
a {\it $3$-partitioned} surface-link.
The definition of the triple linking
of $F=K_1\cup K_2\cup K_3$ is generalized directly
for 3-partitioned surface-links, and
all results and proofs in Sections~2 and 3
are valid for 3-partitioned surface-links.
Theorem~1.5 is a special case of the following:

\begin{theorem}
If two $3$-partitioned surface-links
$F=K_1\cup K_2\cup K_3$ and
$F'=K'_1\cup K'_2\cup K'_3$ are link homologous,
then $\tlk(K_i,K_j,K_k)=\tlk(K'_i,K'_j,K'_k)$.
\end{theorem}

\noindent
{\it Proof.\/}
It is sufficient to prove $\tlk(K_1, K_2, K_3) =
\tlk(K'_1, K'_2, K'_3)$ in a special case that
$K_i=K'_i$, $K_j= K'_j$ and $K_k$ is homologous to
$K'_k$ in $S^4 \backslash (K_i \cup K_j)$,
where $\{i,j,k\}=\{1,2,3\}$.
If $k=2$, then
$\tlk_4(K_1, K_2, K_3)= \tlk_4(K_1, K'_2, K_3)$
by definition.
If $k=1$, then
$\tlk_3(K_1, K_2, K_3)= \tlk_3(K'_1, K_2, K_3)$.
(This is seen as follows: Let $M_1$ be a 3-chain with
$\partial M_1=K_1$.
Since $K'_1$ is homologous to
$K_1$ in $S^4 \backslash (K_2\cup K_3)$, there is a 3-chain
$B$ in $S^4 \backslash (K_2\cup K_3)$
with $\partial B = K'_1 - K_1$.
Let $M'_1=
M_1 +B$, which is a 3-chain with $\partial M'_1= K'_1$.
Then $[M_1 \cap K_2]_{E_3}= [M'_1 \cap K_2]_{E_3}$
in $H_1(E_3) \cong H_1(S^4 \backslash K_3)$.
)
The case $k=3$ is reduced to the previous case ($k=1$)
by use of
Theorem~1.2(i).
\ \qed


\section{Proof of Proposition 1.6}

(1) Let $\ell=k_1\cup k_2$ be
a $(2,2a)$-torus link in a $3$-disk $D^3$
with $\lk(k_1,k_2)=a$.
Let $\gamma$ be a simple loop in ${\bf R}^4$
which intersects a $3$-disk $B_0$ in ${\bf R}^4$
transversely at a single interior point of $B_0$
in the positive direction.
Identify $D^3\times S^1$
with a regular neighborhood $N(\gamma)$
of $\gamma$ in ${\bf R}^4$ and
let $T_1\cup T_2$ be the image of
$\ell\times S^1=k_1\times S^1\cup k_2\times S^1$
in ${\bf R}^4$.
Let $F=K_1\cup K_2\cup K_3$
be a surface-link with
$K_1=T_1$, $K_2=T_2$ and $K_3=\partial B_0$.
Then $F$ is
the desired link.

(2)
Let $\ell= k_1 \cup k_2 \cup k_3$ be a
pretzel link of type $(2a, -2b)$ in a 3-disk $D^3$
so that
$\lk(k_1, k_2)=a$, $\lk(k_2, k_3)= -b$ and
$\lk(k_1, k_3)=0$.
Let $B_1, B_2, B_3$ be mutually disjoint 3-disks
embedded in ${\bf R}^4$
and let $\gamma$ be a simple loop in ${\bf R}^4$
which intersects $B_i$ ($i=1,2,3$) transversely
at
a single interior point of $B_i$ in the positive direction.
Identify $D^3 \times S^1$ with a regular
neighborhood $N(\gamma)$ of $\gamma$ in ${\bf R}^4$ and let
$T_1 \cup T_2 \cup T_3$ be
the image of $\ell \times S^1 =
k_1 \times S^1 \cup k_2 \times S^1 \cup k_3 \times S^1$
in ${\bf R}^4$.
Let $F = K_1 \cup K_2 \cup K_3$ be a surface-link obtained from
$(T_1 \cup T_2 \cup T_3) \cup (\partial B_1 \cup
\partial B_2 \cup \partial B_3)$ by
piping such that $F$ has a projection as in Figure~3.
Then $F$ is
the desired link.
\ \qed




\begin{figure}[htb]
\begin{center}
\mbox{
\epsfxsize=5in
\epsfbox{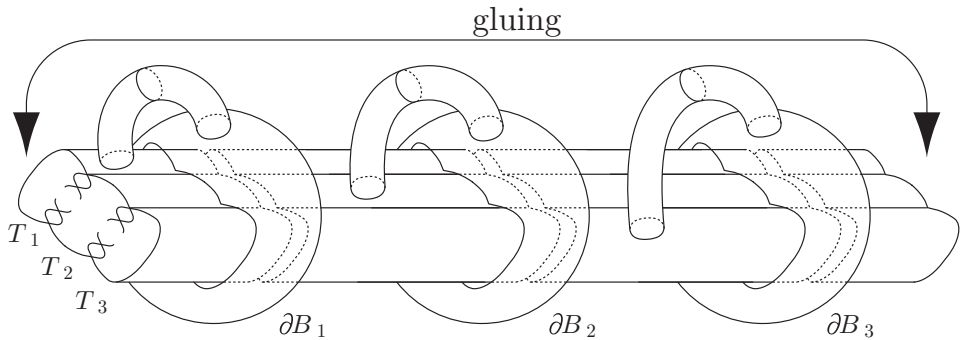}
}
\end{center}
\caption{Linked tori with given linking invariants}
\label{Fig03}
\end{figure}


\section{Remarks}

\noindent
{\bf Remark 6.1} \
The definition of $\tlk_6$
can be seen as a direct analogue of (6) given
in \cite{Rolfsen} page 133.
See also \cite{FennRolf,MassayRolf}.
This generalizes the triple linking to all link maps, instead of
embeddings.
Moreover, it easily generalized to all
dimensions. Let $M_i$ denote a closed connected $n$-manifold for
$i=1,\ldots , n+1$. Let an embedding $f: \cup_{i=1}^{n+1} M_i
\rightarrow  {\bf R}^{n+2}$ be given.
Define $L: \prod_{i=1}^{n+1} M_i \rightarrow
\prod_{j=1}^n S^{n+1}_j$ as
follows. Let $x_i \in M_i$; for $i=1, \ldots, n$, let
$\Delta_i = f(x_i)-f(x_{i+1})/ || f(x_i)-f(x_{i+1}) ||$.
Then
$$L(x_1,\ldots x_{n+1}) = ( \Delta_1, \Delta_2, \ldots , \Delta_n).$$
The general $(n+1)$-fold linking number, $\glk$,
is defined by
$$\glk(f(M_1), \dots, f(M_{n+1}))= \deg(L).$$
We can generalize the notion of $\tlk_1$ to
a diagram in ${\bf R}^{n+1}$ of an
$(n+1)$-component $n$-manifold-link
$M_1 \cup \dots \cup M_{n+1}$
in ${\bf R}^{n+2}$; namely, a diagram has
generic $(n+1)$-tuple points and we count the number
of times $M_1$ is over $M_2$ is over ... is over $M_{n+1}$
with signs.
It is difficult to show that
this value is an invariant of the
$n$-manifold-link in ${\bf R}^{n+2}$ directly,
since we do not know
Reidemeister moves
for higher dimensions $(n \geq 3)$.
However, the proof of Lemma 3.5 goes through to
show that $\glk$ is the same as this count (up to sign).
Thus we have that this number (generalization of $\tlk_1$)
is an invariant
of 
an
$(n+1)$-component $n$-manifold-link.

\bigskip

\noindent
{\bf Remark 6.2} \
In classical link theory,
the linking number determines
the link homology classes completely.
However,
the triple linking of surface-links
is {\it not} a complete invariant
of the surface-link homology;
there exists a pair of surface-links
with the same triple linking invariants
which are not link homologous.
A classification of surface-link homology classes
is discussed in a forthcoming paper.


\bigskip

\noindent
{\large\bf Acknowledgments} \
The authors would like to thank
Akio Kawauchi for helpful suggestions.
JSC is being supported by NSF grant DMS-9988107.
SK and SS are being supported by
Fellowships
from the Japan Society for the Promotion of Science.
MS is being supported by NSF grant DMS-9988101.


\end{document}